\documentclass{amsproc}
\usepackage{amssymb,amsmath,eepic,verbatim}

\newtheorem{theorem}{Theorem}
\newtheorem{corollary}[theorem]{Corollary}

\newtheorem{proposition}[theorem]{Proposition}

\newtheorem{lemma}[theorem]{Lemma}
\newtheorem{lemdef}[theorem]{Lemma/Definition}
\theoremstyle{definition}

\theoremstyle{remark}

\newtheorem{example}{Example}
\renewcommand{\O}{\mathcal{O}}

\newcommand{\C}{\mathbb{C}}

\newcommand{\Z}{\mathbb{Z}}
\newcommand{\Q}{\mathbb{Q}}
\renewcommand{\P}{\mathbb{P}}

\newcommand\lie[1]{\mathfrak{#1}}

\newcommand{\g}{\lie{g}}

\newcommand{\p}{\lie{p}}

\newcommand{\Alc}{\lie{A}}
\renewcommand{\u}{\lie{u}}

\newcommand{\on}{\operatorname}

\newcommand{\Hom}{ \on{Hom}}

\newcommand{\ssm}{\kern-.5ex \smallsetminus \kern-.5ex}

\newcommand\dirac{/\kern-1.2ex\partial} 
\newcommand\qu{/\kern-.7ex/} 
 
\renewcommand{\comment}[1]   {{}}
\newcommand{\labell}\label

\newcommand{\ol}{\overline}

\newcommand{\ti}{\tilde}

\newcommand\E{\mathcal{E}}

\newcommand\Flag{\on{Flag}}

\newcommand\aff{\on{aff}}
\renewcommand{\ss}{\on{ss}}

\newcommand\Gr{\on{Gr}}

\newcommand\rank{\on{rank}}

\begin{document}

\def\mathunderaccent#1{\let\theaccent#1\mathpalette\putaccentunder}
\def\putaccentunder#1#2{\oalign{$#1#2$\crcr\hidewidth \vbox
to.2ex{\hbox{$#1\theaccent{}$}\vss}\hidewidth}}
\def\ttilde#1{\tilde{\tilde{#1}}}

\title[Comparison formula for Gromov-Witten invariants]{On D. Peterson's comparison formula \\ for Gromov-Witten
invariants of $G/P$}

\author{Christopher T. Woodward} \thanks{This research was partially
supported by NSF grants DMS9971357 and DMS0093647.}
\address{Department of Mathematics, Hill Center, Rutgers University, 110
Frelinghuysen Road, Piscataway, New Jersey 08854-8019.}
\email{ctw@math.rutgers.edu}

\begin{abstract} We prove a formula of Dale Peterson comparing
Gromov-Witten (GW) invariants of $G/P$ to those of $G/B$ using
canonical reductions of bundles.
\end{abstract}



\maketitle

An unpublished formula of Dale Peterson describes how $3$-point, genus
$0$ Gromov-Witten invariants of $G/P$ compare with those of $G/B$.
Our purpose in this note is to describe an explanation, and in
particular a proof, of this formula using ideas from moduli of
principal bundles over curves.  The quantum product with respect to
the Schubert basis in $G/B$ can be computed either recursively using
Peterson's quantum Chevalley formula, proved in \cite{fu:qp}, or using
polynomial representatives for the Schubert classes in the
Givental-Kim presentation of the small quantum cohomology
\cite{fo:qu}, \cite{ma:qgb}.  Together these results give a
practicable method for computing the small quantum cohomology in the
Schubert basis for arbitrary $G/P$, although there are much more
effective methods in many special cases
\cite{be:qs,be:qm,cf:qc3,kr:or,kr:la,po:af}.

The idea of the proof is the following.  Given a morphism $\varphi$ of
$\P^1$ to a partial flag variety $X$ of a certain degree $d$, we can
pull back the tautological bundles over $X$.  Giving a lift $\varphi'$
of degree $d'$ of $\varphi$ to a partial flag variety $X'$ dominating
$X$ is equivalent to giving filtrations of the pull-back of the
tautological bundles, by sub-bundles of ranks and degrees determined
by the data $X',d'$.  It turns out that for general $\varphi$ one can
determine the degree $d'$ of the lift corresponding to the
Harder-Narasimhan filtration.  This produces a birational equivalence
between the space of morphisms $\Hom_d(\P^1,X)$ of degree $d$ to $X$,
and the space of morphisms of degree $d'$ to $X'$.  Playing a similar
game with the Jordan-H\"older filtration relates this moduli space to
a moduli space of morphisms of $\P^1$ to the full flag variety.  The
idea for arbitrary $G/P$ is the same but uses the parabolic reductions
of Atiyah-Bott and Ramanathan for principal bundles over curves, which
generalize the Harder-Narasimhan and Jordan-H\"older filtrations for
vector bundles.

We adopt the notation of our joint paper with W. Fulton \cite{fu:qp}.
In particular, $G$ is a connected, simply connected, semisimple
complex Lie group with Borel subgroup $B$, opposite Borel subgroup
$B\overline{\phantom{a}}$, maximal torus $T$, and Weyl group $W$.  Let
$w_o$ be the longest element of $W$.  Let $P$ be a standard parabolic
subgroup, corresponding to a subset $\Delta_P$ of the simple roots.
Let $R_P^+$ denote the set of roots that are combinations of elements
of $\Delta_P$.  For any $u \in W/W_P$, the opposite Schubert variety
is $Y(u) = \overline{B\overline{\phantom{a}}uP/P}.$ Its class in the
integral cohomology ring $H^{\bullet}(G/P)$ is denoted by $\sigma_u$.
The dual cohomology class is $\sigma^u := \sigma_{w_o u}$.  Let $n \ge
3$ be an integer, $p_1,\ldots,p_n \in \P^1$ distinct points, and
$g_1,\ldots,g_n \in G$ general elements.  For any $u_1,\ldots,u_n \in
W/W_P$, define
$$\langle \sigma_{u_1}, \ldots, \sigma_{u_n}
\rangle_d = \# \{ \varphi: \ \P^1
\to G/P, \ \deg(\varphi) = d, \ \varphi(p_i) \in g_i Y(u_i) \ 
\text{for} \ i = 1,\ldots, n \} $$ 
if this number is finite, and zero otherwise.  These invariants may
also be defined as pairings in the Kontsevich-Manin moduli space
$\ol{M}_{0,n}(G/P,d)$ of degree $d$ $n$-pointed genus $0$ stable maps.
Namely, let
$$ f: \ \ol{M}_{0,n}(G/P,d) \to \ol{M}_{0,n}, \ \ e_i:
\ol{M}_{0,n}(G/P,d) \to G/P$$
denote the forgetful morphism to the moduli space of stable
$n$-pointed genus $0$ curves, resp. the $i$-th evaluation map.  Then
$\langle \sigma_{u_1}, \ldots, \sigma_{u_n} \rangle_d$ is the
coefficient of the point class in $f_*( e_1^* \sigma_{u_1} \cdot
\ldots \cdot e_n^* \sigma_{u_n})$.

Define a deformation of the cohomology ring of $G/P$ as follows.  Let
$s_1,\ldots, s_r$ be the simple reflections in $W$ not in $W_P$.  The
classes $\sigma^{[s_1]},\ldots,\sigma^{[s_r]}$ form a basis for
$H^{\dim(G/P)-2}(G/P)$ which we identify with $H_2(G/P)$.  For any
degree $ d = \sum_{i=1}^r d_i \,\sigma^{[s_i]} $ set $q^d = q_1^{d_1}
\cdot \ldots \cdot q_r^{d_r}$ in $\Z[q] := \Z[q_1,\ldots,q_r]$.  The
quantum multiplication formula
$$ \sigma_{u_1} \star \ldots \star \sigma_{u_{n-1}} = \sum_d q^d
\sum_{u_n} \langle \sigma_{u_1}, \ldots, \sigma_{u_n} \rangle_d
\sigma^{u_n}
$$
defines an associative, commutative, $\Z[q]$-linear product on
$$QH^{\bullet}(G/P) =
H^{\bullet}(G/P)\otimes_{\mathbb{Z}}{\mathbb{Z}}[q]$$ 
the {\em small quantum cohomology} ring of $G/P$.  We call the
structure coefficients $\langle \sigma_{u_1}, \ldots, \sigma_{u_n}
\rangle_d$ the small GW-invariants of $G/P$.  These invariants, are
defined using the pull-back of the point class on the moduli space of
stable maps, and should not be confused with the $n$-point
GW-invariants of $G/P$ that play a role in the large quantum
cohomology and are less well understood.

Actually it is somewhat misleading to call the ring
$QH^{\bullet}(G/P)$ cohomology, since it is not functorial: A morphism
$h: \ X \to X'$ does not induce a morphism $QH^{\bullet}(X') \to
QH^{\bullet}(X)$ unless $h$ is an isomorphism.  In particular, the
projection $G/B \to G/P$ does not induce a morphism $QH^{\bullet}(G/P)
\to QH^{\bullet}(G/B)$.  Peterson's comparison formula \eqref{comp}
below fills this gap: it expresses the degree $d_P$ invariants of
$G/P$ in terms of degree $d_B$ invariants for $G/B$.  Unfortunately
the definition of $d_B$, which follows, is not very explicit.  Let
$$\phi_{P/B}: \ G/B \to G/P$$ 
be the projection.  For any weight $\mu$, let $L(\mu)$ denote the
corresponding line bundle over $G/B$ and $c_1(L(\mu)) \in H^2(G/B)$
its first Chern class.  We denote by $( \ , \ )$ the pairing of
homology and cohomology.

\begin{lemdef} \label{deg}  For any $d_P \in H_2(G/P)$, there exists a 
unique $d_B \in H_2(G/B)$ such that $(\phi_{P/B})_* d_B = d_P$ and
$$(d_B, c_1(L(\alpha)) \in \{ 0, 1 \}, \ \ \ \forall \alpha \in R_P^+
.$$
Furthermore, if $\Hom_{d_P}(\P^1,G/P)$ is non-empty then so is
$\Hom_{d_B}(\P^1,G/B)$.
\end{lemdef}

\begin{proof}   Denote by $\pi_B^*$ the
isomorphism from $H^2(G/B)$ to the weight lattice
$$\pi_B^*: \ H^2(G/B) \to \Lambda^*, \ \ \ \ c_1(L(\mu)) \mapsto
\mu $$
and by $\pi_B$ the dual isomorphism $\pi_B: \Lambda \to H_2(G/B)$.
For any parabolic subgroup $P \subset G$ we have similar isomorphisms
$$ \pi^*_P : \ H^2(G/P) \to (\Lambda^*)^{W_P}, \ \ \ \pi_P: \
\Lambda^P \to H_2(G/P) $$
where $\Lambda^P := ((\Lambda^*)^{W_P})^* .$ Let $r_P: \ \Lambda \cong
\Lambda^{**} \to \Lambda^P $ denote the map given by restriction.  Let
$\Lambda_P$ denote the coweight lattice for the semi-simple part of
the Levi factor of $P$, and $W_P^{\aff} = W_P \ltimes \Lambda_P$ the
affine Weyl group for $P$. The inverse image $r_P^{-1}(\lambda_P)$ is
invariant under the action of $W_P^{\aff}$, and
$$\Alc_P = \{ \xi \in \Lambda \otimes_\Z \Q, \ 0 \leq \alpha(\xi) \leq
1, \ \forall \alpha \in R_P^+ \} $$
is a fundamental domain for the action of $W_P^{\aff}$; see
e.g. \cite[p. 90]{hu:re}.  So there is a lift $\lambda_B$ of
$\lambda_P$ in $\Alc_P$. Let $d_B = \pi_B(\lambda_B)$.

It follows from e.g. the discussion in \cite{fu:st} that
$\Hom_{d_P}(\P^1,G/P)$ is non-empty if and only if $d_P$ is a
non-negative combination of the classes $\sigma^{[s_i]}$ for $
\alpha_i \in R_P^+$.  Suppose that the latter holds.  By
e.g. localization \cite[Lemma 2.1]{fu:qp}, $\sigma^{s_i} = \pi_B( -
h_i)$, where $h_i$ is the coroot of $\alpha_i$.  Write $\lambda_B =
c_1 h_1 + \ldots c_n h_n$.  We may assume without loss of generality
that the only positive coefficients are $c_1,c_2,\ldots,c_k$, for some
$k \leq j$.  Let $\alpha$ denote the highest root for the parabolic
subgroup defined by this subset.  Then $(\alpha,h_j) \ge 0$ for $j
\leq k$ and $(\alpha,h_j) \le 0$ for $j > k$.  This implies $
(\lambda_B,\alpha) \ge 2$.  We have $c_i = (\lambda_B,\omega_i) =
(\lambda_P,\omega_i) \leq 0$ for $\alpha_i \notin \Delta_P$, where
$\omega_i$ is the corresponding fundamental weight.  Therefore the
simple roots $\alpha_1,\ldots,\alpha_k$ are in $\Delta_P$ and $\alpha
\in R_P^+$ which contradicts the definition of $\lambda_B$.  This
shows that $\lambda_B$ is a non-positive combination of the simple
coroots, so $d_B$ is a non-negative combination of the classes
$\sigma^{s_i}$, so $\Hom_{d_B}(\P^1,G/B)$ is non-empty.\end{proof}
In some cases one can find simple formulas for $d_B$:

\begin{example} \label{prev}  Suppose $G = SL(3)$, $P = P_{\omega_1}$.  Then
$G/P = \P^2$ and $H_2(G/P) \cong \Z$, with generator $\sigma^{[s_1]} =
[\P^1]$.  Let $h_1,h_2 \in \Lambda$ denote the simple coroots.  Given
a degree $d_P = d_1 \sigma^{[s_1]}$, we have $\lambda_P = - d_1
r_P(h_1)$.  The lifts of $\lambda_P$ are of the form $\lambda_B = -
d_1 h_1 - d_2 h_2$.  To find $d_B$, we solve for $d_2$ so that
$$(\alpha_2,\lambda_B) = d_1 - 2d_2 \in \{ 0, 1 \} .$$
The solution is $d_2 = d_1/2$, if $d_1$ is even, and $d_2 = (d_1
-1)/2$, if $d_1$ is odd.
\end{example} 

Define $P'$ to be the parabolic subgroup of $G$ so that $\Delta_{P'} =
\{ \alpha \in \Delta_P, \ \alpha(\lambda_B) = 0 \}$.  Let $d_{P'}$
denote the image of $d_B$ under the projection $H_2(G/B) \to
H_2(G/P')$ and $\lambda_{P'} = \pi_{P'}^{-1}(d_{P'})$.  Let $w_{P'}$
denote the longest element of the Weyl group $W_{P'}$.  For any $u \in
W/W_{P}$, let $\ti{u} \in W$ denote its minimal length lift.

\begin{theorem}[Peterson's Comparison Formula] \label{compare}
Let $u_1,\ldots,u_n \in W/W_P$.  For any degree $d_P\in H_2(G/P)$
we have for the degree $d_B$ defined by Lemma \ref{deg},
\begin{equation} \label{comp} \langle \sigma_{u_1}, \ldots,  \sigma_{u_n} \rangle_{d_P} = \langle
\sigma_{\ti{u}_1}, \ldots, \sigma_{\ti{u}_{n-1}}, \sigma_{\tilde{u}_n
w_{P'}} \rangle_{d_B} .\end{equation}
\end{theorem}

\begin{example}  Let $G/P = SL(3)/P_{\omega_1} = \P^2$ and
$d_P = \sigma^{[s_1]}$ be the generator of $H_2(\P^2)$.  Then
$\sigma_{[s_1]}$ is the cohomology class of a line and
$\sigma_{[s_2s_1]}$ is the class of a point.  Since there is a unique
line passing through a line and two points in general position in
$\P^2$,
$ \langle \sigma_{[s_1]}, \sigma_{[s_2s_1]}, \sigma_{[s_2s_1]}
\rangle_{d_P} = 1 .$
The lift $d_B = \sigma^{s_1}$ in $H_2(G/B)$, by Example \ref{prev}.
Hence $P' = B$ and $w_{P'} = e$ is the identity in $W$.  One can check
that 
$ \langle \sigma_{s_1}, \sigma_{s_2 s_1}, \sigma_{s_2s_1} \rangle_{d_B} =
1 $
using the Peterson's quantum Chevalley formula \cite{fu:qp}, or
explicitly as follows: The intersection $e_1^{-1}(Y(s_2s_1)) \cap
e_2^{-1}(w_o Y(s_2s_1)) \subset \ol{M}_{0,3}(G/B,d_B)$ is proper, and
maps isomorphically under $e_3$ onto $s_1 Y(s_1s_2)$.  The latter
meets $Y(s_1)$ properly at $x(s_2s_1) \in G/P$, which implies that the
GW-invariant is $1$.  Here $x(s_2s_1)$ denotes the $T$-fixed point
corresponding to $s_2s_1 \in W$.
\end{example} 

We prove Theorem \ref{compare} at the end of the paper using Theorem
\ref{immerse} below.  Recall that the set $\Hom_{d_P}(\P^1,G/P)$ of
degree $d_P$ morphisms $\P^1 \to G/P$ has the structure of a smooth,
quasi-projective variety.  Denote by $\phi_{P'/B}, \ \phi_{P/P'}$ the
projections
\begin{equation} \label{phis}
 \phi_{P'/B}: \ G/B \to G/P', \ \phi_{P/P'}: G/P' \to G/P .\end{equation}
We denote by $\Hom_{d_{P'}}(G/P') \times_{G/P'} G/B$ the fiber product
over $G/P'$ via evaluation at $0$ and $\phi_{P'/B}$.

\begin{theorem} \label{immerse} The morphism
\begin{equation} \label{BP1}
\Hom_{d_B}(G/B) \to \Hom_{d_{P'}}(G/P') \times_{G/P'} G/B, \ \ \varphi
\mapsto (\phi_{P'/B} \circ \varphi, \varphi(0))
\end{equation}
is an open, dense immersion.  The morphism 
\begin{equation} \label{P1P}
 \Hom_{d_{P'}}(G/P') \to \Hom_{d_P}(G/P), \ \ \varphi \mapsto \phi_{P/P'}
\circ \varphi \end{equation}
is birational. 
\end{theorem}

Theorems \ref{compare} and \ref{immerse} were both stated in
\cite{pe:mi} without proof.  We will prove them using basic facts on
semistability of principal bundles over curves.  Recall that a vector
bundle $E \to X$ over a curve $C$ is {\em semistable} if every
sub-bundle $E' \subset E$ has {\em slope} $\mu(E') =
\deg(E')/\rank(E')$ at most the slope $\mu(E)$ of $E$.  If $E$ is not
semistable, there is a unique sub-bundle $E'$ of maximal slope that is
maximal rank among sub-bundles of slope $\mu(E')$.  Applying this fact
inductively leads to the Harder-Narasimhan filtration, which is the
unique filtration with the given degrees and ranks.

In order to make what follows more readable, we will first prove the
theorem for a simple example.  Consider the case that $G = SL(3)$, $P
= P_{\omega_1}$, and $\lambda_P = r_P(h_1)$ so that $d_P$ is the
degree of a line in $G/P = \P^2$.  Over $\P^2$ we have the quotient
vector bundle $Q$ and the tautological bundle $R$, of ranks $2$, $1$
respectively, given by
$$R_{[z]} = [z], \ \ \ \ Q_{[z]} = \C^3/[z] , \ \ \ \ [z] \in \P^2
.$$
Any morphism $\varphi_P: \ \P^1 \to \P^2$ of degree $d_P$ maps $\P^1$
isomorphically onto a line in $\P^2$.  A theorem of Grothendieck
states that any vector bundle splits over $\P^1$; in this example $
\varphi_P^* Q \cong \O(1) \oplus \O(0), \ \ \varphi_P^*R \cong
\O(-1).$ One way of seeing this is to note that $Q \otimes R^{-1}$ is
the tangent bundle $T\P^2$ of $\P^2$; the pull-back $\varphi^* T\P^2$
is the sum of the tangent bundle $T\P^1 \cong \O(2)$ to $\P^1$ and the
normal bundle $N\P^1 \cong \O(1)$.  It follows that the
Harder-Narasimhan filtration of $\varphi^*_P Q$ is has a single
non-trivial term given by the line bundle $ S$ isomorphic to $ \O(1)
.$ The choice of a line sub-bundle of $Q$ defines a lift $\varphi_B: \
\P^1 \to G/B = \Flag(\C^3)$ of $\varphi_P$ as follows.  Let $\pi_{[w]}
: \C^3 \to \C^3/\varphi_P([w])$ denote the projection.  Define
$$ \varphi_B([w]) = \left( \varphi^*_P R_{[w]} \subset \pi^{-1}_{[w]}
S_{[w]} \right) .$$
A little yoga with the definition of degree shows that the element
$\lambda_B = \pi_B^{-1}(\deg(\varphi_B))$ satisfies
$$ ( \lambda_B, \omega_1) = c_1(\varphi_P^* R) = - 1, \ \ \ \ \
(\lambda_B,\omega_2) = c_1(\varphi_P^* R \oplus S) = 0 $$
which implies $\lambda_B = - h_1$.  The fact that the
Harder-Narasimhan filtration is the unique filtration with given
degrees implies that $\varphi_B$ is the unique lift of $\varphi_P$ of
degree $d_B = \pi_B(\lambda_B)$.  Since this is true for any map
$\varphi_B: \P^1 \to \P^2$ of degree $d_P$, the map
$$ \Hom(\P^1,\Flag(\C^3))_{d_B} \to \Hom(\P^1,\P^2)_{d_P}, \ \ \ \
\varphi_B \mapsto \varphi_P := \phi_{P/B} \circ \varphi_B $$
is a bijection.  Since both varieties are smooth it is an isomorphism;
this is a special case of Theorem \ref{immerse}.  In this example, $P'
= B$, so \eqref{BP1} is a tautology.  In general, the proof of
\eqref{BP1} involves the Jordan-H\"older filtration, as we explain
below.

In order to prove Theorem \ref{immerse} in general, we need some
terminology for principal $G$-bundles over a variety $X$.  First, a
principal $G$-bundle $\E \to X$ is a right $G$-variety over $X$ that
is locally trivial; in our situation we may assume local triviality in
the Zariski topology.  For any principal $G$-bundle $\E \to X$ and
morphism $\varphi: \ X' \to X$, we denote by $\varphi^*\E$ the
pull-back bundle.  For any left $G$-variety $F$ we denote the
associated fiber bundle by $\E(F)$.  Let $G' \subset G$ be a subgroup.
A {\em reduction} of $\E$ to $G'$ is a section $\sigma$ of the fiber
bundle $\E(G/G')$.  A special role is played by reductions to maximal
parabolic subgroups $P \subset G$.  In the case $G = GL(V)$, the
maximal parabolic subgroups are the stabilizers of subspaces $V'
\subset V$.  A parabolic reduction $\sigma: \ X \to \E(G/P)$ is
equivalent to a sub-bundle of the associated vector bundle $\E(V)$
with fiber $V'$.

Semistability of principal $G$-bundles is defined as follows.  For any
standard maximal parabolic $P$, let $\omega_P$ be the fundamental
weight such that $\Delta_P$ is the set of simple roots vanishing on
$\omega_P$.  A principal $G$-bundle $\E \to X$ is called semistable if
and only if for any reduction $\sigma: \ C \to \E/P$ to a standard
maximal parabolic $P$, the degree of the associated line bundle
$\sigma^* \E(\omega_P)$ is non-positive.  For $G = SL(n)$,
semistability of $\E$ is equivalent to semistability of the associated
vector bundle (see Ramanathan \cite{ra:th} or Atiyah-Bott
\cite[Section 10]{at:mo}). For any $G$, semistability of $\E$ is
equivalent to semistability of the vector bundle $\E(\g)$ associated
to the adjoint representation $\g$.  If $\E$ is not semistable, there
is a canonical {\em Atiyah-Bott} parabolic reduction $\sigma_{\E}: \ C
\to \E/P_\E$, where the parabolic subgroup $P_\E$ has Lie algebra
$\p_\E$ isomorphic to the fiber of the degree-zero term $\E(\g)_0$ in
the Harder-Narasimhan filtration of $\E(\g)$.  The canonical reduction
has a uniqueness property generalizing that of the Harder-Narasimhan
filtration: For any reduction $\sigma: \ X \to \E/P$, define the slope
$\mu$ of $\sigma$ to be the homomorphism from characters $\chi$ of $P$
to $\Z$ given by mapping $\chi$ to the degree of the associated line
bundle $\sigma_\E^* \E(\chi)$.
\begin{proposition} (see e.g. \cite[pp.11-12]{tw:pa})  \label{uniq}
$\sigma_\E$ is the unique reduction of $\E$ to $P_\E$ with slope
$\mu_\E$.
\end{proposition}

If a degree $0$ vector bundle $E \to C$ is semistable, there is a
Jordan-H\"older filtration on $E$ characterized by the property that
the associated graded bundle $\Gr(E)$ is semistable, and the
filtration is maximal among filtrations of this type.  The
Jordan-H\"older filtration is not unique; however, $\Gr(E)$ is unique
up to isomorphism.  The corresponding notion for principal bundles was
introduced by Ramanathan \cite{ra:th}: A reduction $\sigma: \ C \to
\E/P$ is called {\em admissible} if $\sigma$ has slope $0$.  Let $L$
denote the standard Levi subgroup of $P$ and $\pi_L: P \to L$ and
$\iota_L: \ L \to G$ denote the homomorphisms given by projection and
inclusion respectively.

\begin{proposition} \label{admiss}  
\cite[3.5.11]{ra:th} Let $\sigma: \ C \to \E/P$ be an admissible
reduction of $\E$.  $(\iota_L)_* (\pi_L)_* \sigma^*\E$ is semistable
if and only if $\E$ is semistable.
\end{proposition}
\noindent If $\sigma$ is admissible and $(\pi_L)_* \sigma^*\E$ is
stable, call $\sigma$ a {\em Ramanathan} reduction.  By
\cite[Proposition 3.12]{ra:th}, Ramanathan reductions exist for any
bundle $\E$.  Define an equivalence relation on principal $G$-bundles
by $\E \sim (\iota_L)_* (\pi_L)_* \sigma^* \E(G)$, where $\sigma$ is a
Ramanathan reduction.  Ramanathan \cite{ra:th} constructs a coarse
moduli space for equivalence classes of semistable principal bundles.
In genus zero, the moduli problem is trivial, for the following reason
which is an easy consequence of Grothendieck's theorem that any
principal $G$-bundle over $\P^1$ admits a reduction to $T$
\cite{gr:cl}:
\begin{theorem}  \label{BG} Any semistable principal $G$-bundle
$\E \to \P^1$ is trivial: $\E \cong \P^1 \times G$.
\end{theorem} 
\noindent These results have straightforward generalizations to
the case that $G$ is reductive.

We apply these results to pull-backs of bundles on $G/P$.  Let
$\varphi_P: \ X \to G/P$ be a morphism and $\E_P$ the principal
$P$-bundle $G \to G/P$.  For any parabolic subgroup $P' \subseteq P$,
lifts $\varphi_{P'} : \ X \to G/P'$ of $\varphi_P: \ X \to G/P$ are in
one-to-one correspondence with reductions $\sigma_{P'}: \ X \to
\varphi_P^* \E_P(P/P')$.

Our goal is to prove Theorem \ref{immerse} by thinking of it as a
statement about reductions of bundles.  Let $P = LU$ and $P' = L'U'$
denote the standard Levi decompositions.  We study the semistability
of the principal $L$-bundle $(\pi_L)_* \varphi_P^* \E_P$.
\begin{lemma} \label{AB} Suppose there exists a lift $\varphi_B: \ 
\P^1 \to G/B$ of $\varphi_P$ of degree $d_B$.  Then the Atiyah-Bott
canonical reduction of $(\pi_L)_* \varphi_P^* \E_P$ corresponds to the
lift $\varphi_{P'} = \phi_{P'/B} \circ \varphi_B $.
\end{lemma}

\begin{proof}  
Let $B$ act on $L$ via $\pi_L$.  Because of the isomorphisms
$$ \varphi_P^*\E_P(P/P') \to \varphi_B^* \E_B(P/P') \to
\varphi_B^*\E_B(L/L \cap P'), $$
the map $\varphi_{P'}$ defines a reduction of $\varphi_B^*\E_B(L)$ to
$L \cap P'$.  The filtration $\lie{l} \cap \u' \subset \lie{l} \cap
\p' \subset \lie{l}$ is $B$-stable.  We claim that
\begin{equation} \label{filter} 
 \varphi^*_B \E_{B}(\lie{l} \cap \u') \subseteq \varphi^*_B
\E_{B}(\lie{l} \cap \p') \subseteq \varphi^*_B \E_{B}(\lie{l})
\end{equation}
is the Harder-Narasimhan filtration of $\varphi^*_B \E_{B}(\lie{l})$.
We have $ \deg \varphi^*_B \E_{B}(\lie{l}_\mu) = (\lambda_B,\mu). $
Using the definition of the Peterson lift, if $\mu$ is a positive
(resp. negative) root of $\lie{l}$ that is not a root of $\lie{l}'$
then $(\lambda_B,\mu) = 1$ resp. $-1$; otherwise $(\lambda_B,\mu) =
0$.  It follows that the Harder-Narasimhan filtration is
\eqref{filter}, and has slope-zero term $\varphi^*_B \E_{B}(\lie{l}
\cap \p')$.
\end{proof}

\begin{corollary}  Suppose that $\varphi_P$ lifts to a map 
$\varphi_B: \ \P^1 \to G/B$ of degree $d_B$.  Then the composition
$\varphi_{P'}$ of $\varphi_B$ with the projection to $G/P'$ is the
unique lift of $\varphi_P$ to $G/P'$ of degree $d_{P'}$.
\end{corollary}

\begin{proof} By Lemma \ref{AB} and Proposition \ref{uniq}.
\end{proof}

We now consider the comparison between $G/P'$ and $G/B$.  Let
$\varphi_{P'}: \ \P^1 \to G/P'$ be a morphism of degree $d_{P'}$.  Let
$L' \subset P'$ be the standard Levi subgroup of $P'$, $Z(L')$ its
center, and $L'_{\ss} = L'/Z(L')$.  Let $\pi_{L'_{\ss}}: \ P' \to
L'_{\ss}$ denote the projection, and $B'_{\ss}$ the image of $B \cap
L'$ under $\pi_{L'_{\ss}}$. Since both the standard unipotent subgroup
$U' \subset P'$ and $Z(L')$ act trivially on $P'/B$, we have $
\E_{P'}(P'/B) \cong \E_{P'} (L'_{\ss}/B'_{\ss}).$
\begin{lemma} \label{JH}  Suppose that there exists a lift $\varphi_B$
of $\varphi_{P'}$ to $G/B$ of degree $d_B$.  Then the corresponding
reduction $\sigma_B: \ \P^1 \to \varphi_{P'}^*
\E_{P'}(L'_{\ss}/B'_{\ss})$ is a Ramanathan reduction of
$\varphi_{P'}^*\E_{P'}(L'_{\ss})$.
\end{lemma}

\begin{proof}
Any weight for $L'_{\ss}$ defines a weight $\mu$ for $L'$ in the span
of the roots of $L'$.  Hence $(\lambda_B,\mu) = 0$ and the line bundle
$ \varphi_B^* L(\mu) \cong \sigma_B^* \varphi_{P'}^* \E_{P'} (\mu) $
is trivial.  This implies that $\sigma_B$ is admissible.
\end{proof}

\begin{corollary} 
If there exists a lift $\varphi_B$ of $\varphi_{P'}$ to $G/B$ of
degree $d_B$, then the bundle $\varphi_{P'}^* \E_{P'}(P'/B)$ is
trivial.
\end{corollary}

\begin{proof} 
By Lemma \ref{JH}, Theorem \ref{BG} and Proposition \ref{admiss}.
\end{proof}

Now we prove Theorem \ref{immerse}.  The morphism \eqref{BP1} is an
injection.  Indeed, by Lemma \ref{JH} any lift $\varphi_{B}$ gives a
Ramanathan reduction of $\varphi_B^*\E_B(L'_{\ss})$.  A Ramanathan
reduction of the trivial bundle $\P^1 \times L'_{\ss}$ is a constant
morphism $\P^1 \to L'_{\ss}/B'_{\ss}$, and is therefore specified
uniquely by its value at any point in $\P^1$.  The dimension of
$\Hom_{d_B}(G/B)$ is
\begin{eqnarray*}
\dim(\Hom_{d_B}(\P^1,G/B)) &=& \dim(G/B) + (c_1(G/B),d_B) \\
&=& \dim(G/B) + \sum_{\alpha \in R^+} (\alpha,\lambda_B) \\
&=& \dim(G/P') + \dim(P'/B) + \sum_{\alpha \in R^+ \ssm R^+_{P'}}
(\alpha,\lambda_B) \\
&=& \dim(G/P') + (c_1(G/P'),d_{P'}) + \dim(P'/B) \\
&=& \dim(\Hom_{d_{P'}}(\P^1,G/P')) + \dim(P'/B).
\end{eqnarray*}
It follows that \eqref{BP1} is injective.  Since the domain and
codomain are smooth, irreducible (\cite{ki:co},\cite{th:ir}) and the
same dimension, \eqref{BP1} is an open, dense immersion.

Similarly, by Lemma \ref{AB}, the morphism \eqref{P1P} is injective on
the image of \eqref{BP1}.  The domain and codomain have the same
dimension, since
\begin{eqnarray*}
\dim(\Hom_{d_{P'}}(\P^1,G/P')) &=& \dim(G/P') + \sum_{\alpha \in
R^+_{P'}} (\alpha,\lambda_{B}) \\
&=& \dim(G/P') + \sum_{\alpha \in R^+_{P}} (\alpha,\lambda_B) - \#
R^+_P \ssm R^+_{P'} \\
&=& \dim(G/P) + \sum_{\alpha \in R^+_P} (\alpha,\lambda_B) \\ &=&
\dim(\Hom_{d_P}(\P^1,G/P)) .
\end{eqnarray*}
Since the varieties are smooth and irreducible, \eqref{P1P} is an
open, dense immersion on an open subset, and therefore birational.

Theorem \ref{immerse} and Lemma \ref{AB} imply the following curious fact.
\begin{proposition}  \label{curious} 
For general $\varphi_P \in \Hom_{d_P}(\P^1,G/P)$, the pull-back
$\varphi_P^* \E_P(L)$ is semistable if and only if $\lambda_B$ is
$W_P$-fixed, that is, $P = P'$.
\end{proposition}
\begin{example} Let $G = SL(3,\C)$ and $P = P_{\omega_1}$.  Under
the correspondence between principal bundle and vector bundles, the
bundle $\E_P(L)$ corresponds to $Q \oplus R$.  Since (semi)stability
is preserved by tensoring with line bundles, semistability of $Q$ is
equivalent to semistability of $T\P^2$.  Therefore, a general degree
$d_P$ morphism $\varphi_P: \ \P^1 \to \P^2$ has $\varphi^*_P T\P^2$
semistable (hence trivial) if and only if $d_P$ is even.  
\end{example}

Now we prove Theorem \ref{compare}.  Recall the maps $\phi_{P'/B},
\phi_{P/P'}$ from \eqref{phis}.  For any $u \in W/W_P$, we have the
identities
\begin{equation} \label{ident} 
(\phi_{P/B})^* \sigma_{u} = \sigma_{\ti{u}}, \ \ \ \ \ \ \sigma_{u} =
(\phi_{P/B})_* \sigma_{\ti{u}w_{P}}.
\end{equation}
Composing with the projection and collapsing the unstable components
produces morphisms
$$ h_{P'/B}: \ \ol{M}_{0,n+1}(G/B,d_B)
\to \ol{M}_{0,n+1}(G/P',d_{P'}) \times_{G/P'} G/B ,$$
$$ h_{P/P'}: \ \ol{M}_{0,n+1}(G/P',d_{P'}) \to \ol{M}_{0,n+1}(G/P,d_P).
$$
The existence of $h_{P/P'}$ is proved by the same arguments that
construct the forgetful morphism $f$, see \cite{fu:st}.  Theorem
\ref{immerse} implies that these morphisms are birational.  Let
$\phi_1,\phi_2$ denote the projections so that
$$
\phi_1 \times \phi_2: \ \ol{M}_{0,n+1}(G/P',d_{P'}) \times_{G/P'} G/B
\to \ol{M}_{0,n+1}(G/P',d_{P'}) \times G/B $$  
is the canonical inclusion.  Let $u'_j \in W/W_{P'}$ denote the coset
of $\ti{u}_j$.  We denote by superscript $^B$ objects, maps etc. for
$G/B$, and by $^{P'}$ those for $G/P'$.  From \eqref{ident} and the
identities 
$$\phi_{P'/B} \circ e_i^B = e_i^{P'} \circ \phi_1 \circ
h_{P'/B} , \ \ \ \ f^B = f^{P'} \circ \phi_1 \circ h_{P'/B} $$
it follows that for any $w \in W$,
\begin{multline}
f^B_* ( (e_1^B)^* \sigma_{\ti{u}_{1}} \cdot \ldots \cdot (e_{n-1}^B)^*
\sigma_{\ti{u}_{n-1}} \cdot (e_n^B)^* \sigma_w) \\ =
f^{P'}_*( (e_1^{P'})^*
\sigma_{{u}'_{1}} \cdot \ldots \cdot (e_{n-1}^{P'})^* \sigma_{{u}'_{n-1}}
\cdot (e_n^{P'})^* (\phi_{P'/B})_* \sigma_w).
\end{multline}
In particular,
$$
 f^B_* ( (e_1^B)^* \sigma_{\ti{u}_{1}} \cdot \ldots \cdot
(e_{n-1}^B)^* \sigma_{\ti{u}_{n-1}} \cdot (e_n^B)^* \sigma_{\ti{u}_n
w_{P'}}) = f^{P'}_*( (e_1^{P'})^* \sigma_{{u}'_{1}} \cdot \ldots \cdot
(e_{n}^{P'})^* \sigma_{{u}'_{n}}) .$$
Taking the coefficient of the point class in
$H^{\bullet}(\ol{M}_{0,n})$ gives
$$\langle \sigma_{u'_1}, \ldots, \sigma_{u'_n} \rangle_{d_{P'}} =
\langle \sigma_{\ti{u}_1}, \ldots, \sigma_{\ti{u}_{n-1}},
\sigma_{\tilde{u}_n w_{P'}} \rangle_{d_B} .$$ 
A similar but easier
argument shows $ \langle \sigma_{u_1}, \ldots \sigma_{u_n}
\rangle_{d_P} = \langle \sigma_{u'_1}, \ldots, \sigma_{u'_n}
\rangle_{d_{P'}}, $
which completes the proof.

\def\cprime{$'$} \def\cprime{$'$} \def\cprime{$'$}
  \def\polhk#1{\setbox0=\hbox{#1}{\ooalign{\hidewidth
  \lower1.5ex\hbox{`}\hidewidth\crcr\unhbox0}}}



\end{document}